\documentclass[10pt]{article}
\usepackage{amsmath}
\numberwithin{equation}{section}
\usepackage{amssymb}
\usepackage{mathrsfs}
\usepackage{amsfonts}
\usepackage{mathrsfs,amscd,amssymb,amsthm,amsmath,bm,graphicx}
\usepackage{graphicx}
\usepackage{cite}
\usepackage{color}
\usepackage{enumitem}
\usepackage{setspace}
\makeatletter

\renewcommand{\@seccntformat}[1]{{\csname the#1\endcsname}{\normalsize.}\hspace{.5em}}
\makeatother

\def \[{\begin{equation}}
\def \]{\end{equation}}

\newtheorem{thm}{Theorem}[section]

\newtheorem{lem}[thm]{Lemma}

\begin{document}

\title{The (degree)-Kirchhoff index of linear crossed octagonal-quadrilateral networks\footnote
{Supported by the National Natural Science Foundation of China (11771362), Anhui Provincial Natural Science Foundation (2008085J01), and the Natural Science Fund of Education Department of Anhui Province (KJ2020A0478).}
}

\author{Jia-Bao Liu$^{1}$, Ting Zhang$^{1}$, Wenshui Lin$^{2,}$\footnote{Corresponding author. E-mail address: wslin@xmu.edu.cn (W. Lin).}\\
  {\small $^1$ School of Mathematics and Physics, Anhui Jianzhu University, Hefei 230601, P.R. China}\\
  {\small $^2$ School of Informatics, Xiamen University, Xiamen 361005, China}
}

\date{\small (Received May 27, 2022)}
\maketitle

\vspace*{2mm}

\noindent{\bf Abstract}\\
The Kirchhoff index and degree-Kirchhoff index have attracted extensive attentions due to its practical applications in complex networks, physics, and chemistry. In 2019, Liu et al. [Int. J. Quantum Chem. 119 (2019) e25971] derived the formula of the degree-Kirchhoff index of linear octagonal-quadrilateral networks.
In the present paper, we consider linear crossed octagonal-quadrilateral networks $Q_n$. Explicit closed-form formulas of the Kirchhoff index, the degree-Kirchhoff index,
and the number of spanning trees of $Q_n$ are obtained.
Moreover, the Kirchhoff index (resp. degree-Kirchhoff index) of $Q_n$ is shown to be almost 1/4 of its Wiener index (resp. Gutman index).

\vspace*{2mm}

\baselineskip=0.30in

\vspace*{2mm}
\section{Introduction}

Let $G=(V,E)$ be a simple connected graph, where $V = \{v_1, v_2, \ldots, v_n\}$.
The adjacency matrix of $G$ is a 0-1 matrix $A(G)=(a_{ij})_{n \times n}$ with $a_{ij}=1$ iff $v_i v_j \in E$.
Let $D(G)=diag(d_1, d_2, \ldots, d_n)$ be the degree matrix of $G$, where $d_i$ is the degree of vertex $v_i$.
Then $L(G)=D(G)-A(G)$ is termed as the Laplacian matrix,
and $\mathcal{L}(G) =D(G)^{-\frac{1}{2}} L(G) D(G)^{-\frac{1}{2}}$
the normalized Laplacian matrix of graph $G$.
Let $0= \mu_1 <\mu_2 \leq \cdots \leq \mu_n$ be the eigenvalues of $L(G)$,
and $0= \nu_1 <\nu_2 \leq \cdots \leq \nu_n$ the eigenvalues of $\mathcal{L}(G)$.
The sets $Sp(L(G)) = \{\mu_1, \mu_2, \ldots, \mu_n\}$
and $Sp(\mathcal{L}(G)) = \{\nu_1, \nu_2, \ldots, \nu_n\}$ are called
the Laplacian spectrum and normalized Laplacian spectrum of $G$, respectively.

Let $d_{ij}$ denotes the distance between vertices $v_i$ and $v_j$ in $G$, i.e.,
the length of a shortest path connecting them.
The Wiener index \cite{b1} and Gutman index \cite{b2},
defined as $W(G) =\sum_{i<j} d_{ij}$ and $Gut(G) =\sum_{i<j} d_i d_j d_{ij}$,
are two famous distance-based topological indices,
which have successful applications in chemistry and communication networks.

By regarding graph $G$ as an electronic network, in which a unit resistor is placed on each edge of $G$,
Klein and Randi\'{c} \cite{b3} proposed the concept of resistance distance.
The resistance distance $r_{ij}$ between vertex $v_i$ and $v_j$ is defined
as the effective resistance between them.
Similar to the Wiener index and Gutman index, the Kirchhoff index \cite{b3} and
degree-Kirchhoff index \cite{b4} of $G$ are defined as
$Kf(G)=\sum_{i<j} r_{ij}$ and $Kf^*(G)= \sum_{i <j} d_i d_j r_{ij}$, respectively.
Furthermore, $Kf(G)$ (resp. $Kf^*(G)$) was shown to have close relationship
with $Sp(L(G))$ (resp. $Sp(\mathcal{L}(G))$).

\noindent \textbf{Lemma 1.1} \cite{b5,b6}\textbf{.} Let $G$ be a simple graph of order $n \geq 2$.
Then
$$Kf(G) = n \sum_{i=2}^n \frac{1}{\mu_i}.$$

\noindent \textbf{Lemma 1.2} \cite{b4}\textbf{.}
Let $G$ be a simple connected graph of order $n \geq 2$ and size $m$. Then
$$Kf^*(G) = 2m \sum_{i=2}^n \frac{1}{\nu_i}.$$

Due to their practical applications in complex networks, physics, and chemistry, the Kirchhoff index and degree-Kirchhoff index have attracted a lot of attentions.
Up to now, the closed-form formulas of (degree-)Kirchhoff index have been established
for cycles \cite{cycle}, circulant graphs \cite{circulant}, regular graphs \cite{regular}, composite graphs \cite{composite}, complete multipartite graphs \cite{complete-mult}, flower graphs \cite{flower},
the strong prism of a star \cite{strong}, the Cartesian product of $S_n$ and $K_2$ \cite{Cartesian},
zigzag polyhex nanotube TUHC$[2n,2]$ \cite{nano}, linear phenylenes \cite{LinearPheny,LinearPheny2},
cyclic phenylenes \cite{CyclicPheny}, M\"{o}bius phenylenes chain and cylinder phenylenes chain \cite{MobiusPheny,CylinderPheny}, linear [n] phenylenes \cite{Linear-n-Pheny},
generalized phenylenes \cite{GeneralPheny,GeneralPheny2}, periodic linear chains \cite{Periodic},
linear polyomino chain \cite{Linear4}, linear crossed polyomino chain \cite{LinearCrossed4},
linear pentagonal chain \cite{Linear5}, linear hexagonal chain \cite{Linear6,Linear6-2},
linear crossed hexagonal chain \cite{LinearCrossed6},
M\"{o}bius hexagonal chain \cite{Mobius6}, linear octagonal chain \cite{Linear8},
linear crossed octagonal chain \cite{LinearCrossed8}, M\"{o}bius and cylinder octagonal chains \cite{Mobius8}, linear octagonal-quadrilateral chain \cite{Linear8-4},
and M\"{o}bius octagonal-quadrilateral networks \cite{Mobius8-4}.

Let $L_{n}$ be the linear octagonal-quadrilateral network with $n$ quadrilaterals and $n$ octagons, as shown in Figure 1. The crossed octagonal-quadrilateral network $Q_{n}$ (shown in Figure 2) is obtained from $L_n$ by connecting vertex $i$ to vertex $(i+1)'$ and $i'$ to vertex $i+1$, $1 \leq i \leq 4n$.
Inspired by the above works, in this paper we consider the resistance distance-based invariants of linear crossed octagonal-quadrilateral networks.
Explicit closed-form formulas of $Kf(Q_n)$, $Kf^*(Q_n)$, and $\tau(Q_n)$ are obtained, where $\tau(G)$ denotes the number of spanning trees of graph $G$.
Moreover, it is shown that $\lim_{n\rightarrow\infty} Kf(Q_{n})/W(Q_{n})= \lim_{n\rightarrow\infty} Kf^*(Q_{n})/Gut(Q_n) =1/4$.

\begin{figure}[htbp]
\centering\includegraphics[width=15cm,height=4cm]{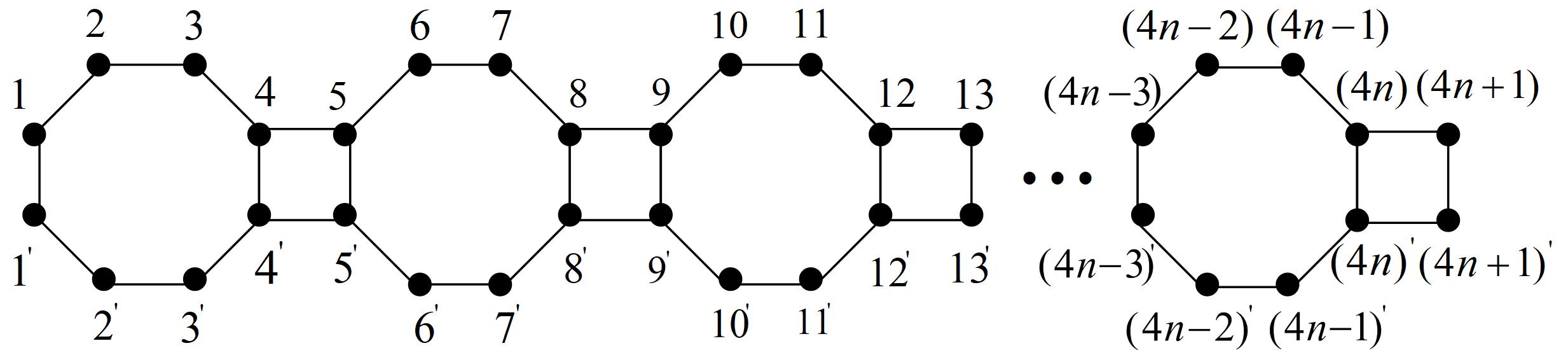}
\caption{The linear octagonal-quadrilateral network $L_n$.}
\end{figure}

\begin{figure}[htbp]
\centering\includegraphics[width=15cm,height=4cm]{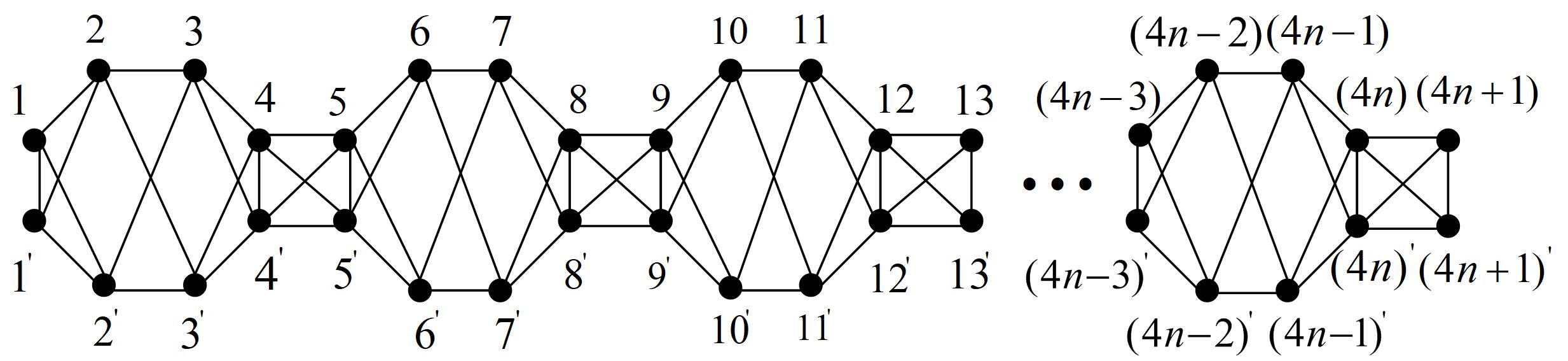}
\caption{The linear crossed octagonal-quadrilateral network $Q_n$.}
\end{figure}

\section{Preliminary}
It is easily seen that $|V(Q_{n})|=8n+2,|E(Q_{n})|=18n+1$, and $\pi=(1,1')(2,2')\cdots\big(4n+1,(4n+1)'\big)$ is an automorphism of $Q_{n}$. Let $ V_{1}=\{1,2,\ldots,4n+1\}$ and $V_2=\{1',2',\ldots,(4n+1)'\}$.
Then we have
\begin{equation*}
  L(G)=\begin{pmatrix}
    L_{V_1V_1} &L_{V_1V_2}\\
    L_{V_2V_1} &L_{V_2V_2}
  \end{pmatrix} ~\text{and}~
  \mathcal{L}(G)=\begin{pmatrix}
    \mathcal L_{V_1V_1} &\mathcal L_{V_1V_2}\\
    \mathcal L_{V_2V_1} &\mathcal L_{V_2V_2}
  \end{pmatrix},
  \end{equation*}
\noindent where $L_{V_i V_j}$ (resp. $\mathcal L_{V_i V_j}$) is the matrix formed by the rows and columns of $L(G)$ (resp. $\mathcal{L}(G)$) corresponding to the vertices in $V_i \bigcup V_j$, $1 \leq i, j \leq 2$.
By symmetry, it holds that $L_{V_1 V_2}=L_{V_2 V_1}$, $L_{V_1 V_1}=L_{V_2 V_2}$, $\mathcal L_{V_1 V_2}=\mathcal L_{V_2 V_1}$, and $\mathcal L_{V_1 V_1}=\mathcal L_{V_2 V_2}$.

Let
\begin{equation*}
  T=\begin{pmatrix}
    \frac{1}{\sqrt{2}}I_{4n+1} &\frac{1}{\sqrt{2}}I_{4n+1} \\
  \frac{1}{\sqrt{2}}I_{4n+1} &-\frac{1}{\sqrt{2}}I_{4n+1}
  \end{pmatrix}.
\end{equation*}

\noindent Then we have

\begin{equation*}
  TL(G)T'=\begin{pmatrix}
    L_A &0\\
    0 &L_S\\
  \end{pmatrix} ~\text{and}~
   T\mathcal{L}(G)T'=\begin{pmatrix}
    \mathcal L_{A} &0\\
    0 &\mathcal L_{S}\\
  \end{pmatrix},
\end{equation*}
\noindent where $L_A= L_{V_{11}}+L_{V_{22}}$, $L_S=L_{V_1V_1}-L_{V_1V_2}$,
$\mathcal L_{A}= \mathcal L_{V_{11}}+\mathcal L_{V_{22}}$, and $\mathcal L_{S}=\mathcal L_{V_1V_1}-\mathcal L_{V_1V_2}$.

Based on the above arguments, by applying the decomposition theorem obtained in \cite{Decom},
we immediately have the following result, where $\Phi_B(x) = det(xI-B)$ refers
the characteristic polynomial of a square matrix $B$.

\begin{lem}
Let $L_{A},L_{S},\mathcal L_{A}, and \mathcal L_{S}$ be defined as above. Then
\begin{eqnarray*}
\Phi_{L(Q_{n})}(x)=\Phi_{L_{A}}(x) \Phi_{L_{S}}(x) ~\text{and}~
\Phi_{\mathcal{L}(Q_{n})}(x)=\Phi_{\mathcal L_{A}}(x)\Phi_{\mathcal L_{S}}(x).
\end{eqnarray*}
\end{lem}

\begin{lem}\textup{\cite{TreeNum}}
Assume that $G$ is a connected graph with $n \geq 2$ vertices,
then $$\tau(G) = \frac{1}{n}\prod_{i=2}^{n} \mu_i.$$
\end{lem}

\section{The Kirchhoff index and Wiener index of $Q_{n}$}\label{sct3}
It is easily seen that
\begin{eqnarray*}
 L_{V_1 V_1}&=&
\left(
  \begin{array}{ccccccccccc}
    3 & -1 & & & & & & & & &\\
    -1 & 4 & -1 & & & & & & & &\\
    & -1 & 4 & -1 & & & & & & &\\
    & & -1 & 5 & -1 & & & & & &\\
    & & & -1 & 5 & -1 & & & & &\\
    & & & & -1 & 4 & -1 & & & &\\
    & & & & & & \ddots & & & &\\
    & & & & & & -1 & 4 & -1 & &\\
    & & & & & & & -1 & 4 & -1 &\\
    & & & & & & & & -1 & 5 & -1\\
    & & & & & & & & & -1 & 3\\
  \end{array}
\right)_{(4n+1)\times(4n+1)}
\end{eqnarray*}
and
\begin{eqnarray*}
 L_{V_1 V_2}&=&
\left(
  \begin{array}{ccccccccccc}
    -1 & -1 & & & & & & & & &\\
    -1 & 0 & -1 & & & & & & & &\\
    & -1 & 0 & -1 & & & & & & &\\
    & & -1 & -1 & -1 & & & & & &\\
    & & & -1 & -1 & -1 & & & & &\\
    & & & & -1 & 0 & -1 & & & &\\
    & & & & & & \ddots & & & &\\
    & & & & & & -1 & 0 & -1 & &\\
    & & & & & & & -1 & 0 & -1 &\\
    & & & & & & & & -1 & -1 & -1\\
    & & & & & & & & & -1 & -1\\
  \end{array}
\right)_{(4n+1)\times(4n+1)}.
\end{eqnarray*}

\noindent Thus,
\begin{eqnarray*}
 L_{A}&=&
\left(
  \begin{array}{ccccccccccc}
    2 & -2 & & & & & & & & &\\
    -2 & 4 & -2 & & & & & & & &\\
    & -2 & 4 & -2 & & & & & & &\\
    & & -2 & 4 & -2 & & & & & &\\
    & & & -2 & 4 & -2 & & & & &\\
    & & & & -2 & 4 & -2 & & & &\\
    & & & & & & \ddots & & & &\\
    & & & & & & -2 & 4 & -2 & &\\
    & & & & & & & -2 & 4 & -2 &\\
    & & & & & & & & -2 & 4 & -2\\
    & & & & & & & & & -2 & 2\\
  \end{array}
\right)_{(4n+1)\times(4n+1)}
\end{eqnarray*}
and $L_S=diag(4,4,4,6,6,4,\ldots,4,4,6,4)$.

Suppose $Sp(L_{A}) = \{0=\rho_{1}, \rho_{2}, \ldots, \rho_{4n+1}\}$
and $Sp(L_S) = \{\sigma_{1}, \sigma_{2}, \ldots, \sigma_{4n+1}\}$. From Lemma 1.1 we have
\begin{equation}
  Kf(Q_{n})=(8n+2)\Big(\sum_{i=2}^{4n+1}\frac{1}{\rho_{i}}+\sum_{i=1}^{4n+1}\frac{1}{\sigma_{i}}\Big).
\end{equation}

Obviously,
\begin{equation}
\sum_{i=1}^{4n+1}\frac{1}{\sigma_{i}}=\frac{1}{4}\times(2n+2)+\frac{1}{6}\times(2n-1)=\frac{5n+2}{6}.
\end{equation}

To get the value of $\sum_{i=2}^{4n+1}\frac{1}{\rho_{i}}$, suppose $\Phi_{L_{A}}(x)=det(xI-L_{A})=x^{4n+1}+a_{1}x^{4n}+\cdots+a_{4n}x$.
Then $1/\rho_2, 1/\rho_3, \ldots, 1/\rho_{4n+1}$ are the roots of the polynomial $a_{4n}x^{4n}+a_{4n-1}x^{4n-1}+\cdots+a_{1}x+1$. From the Vieta theorem we have
\begin{equation}
\sum_{i=2}^{4n+1}\frac{1}{\rho_{i}}=\frac{(-1)^{4n-1}a_{4n-1}}{(-1)^{4n}a_{4n}}.
\end{equation}

Let $W_{i}$ (resp. $U_{i}$) be the $i$-th order principal submatrix of $L_{A}$ formed by its first (resp. last) $i$ rows and $i$ columns, $w_i = det W_i$, and $u_i = det U_i$, $i=1,2,\ldots,4n$.
By direct calculations we have $w_1=2$, $w_2=4$, and $w_i=4w_{i-1}-4w_{i-2}$ if $3\leq i\leq4n$. By solving the recurrence equation we get $w_{i}=2^{i}$, $i=1,2,\ldots,4n$. Furthermore, from the structure of $L_A$, obviously it holds that $u_i = w_i$. Also, let\\
\begin{equation*}
  N_i = \begin{pmatrix}
     4 &-2 & & & & & & & & \\
     -2& 4 & -2 & & & & & & &\\
     &-2 & 4 &-2 & & & & & & \\
     & & -2 &4 & & & & & & \\
     & & & & & \ddots& & & &\\
     & & & & & &  4& -2 & &\\
     & & & & & &  -2& 4 & -2&\\
     & & & & & & &-2 & 4 &-2 \\
     & & & & & & & & -2 & 4
  \end{pmatrix}_{i \times i}
\end{equation*}
\noindent and $n_i=detN_i$, $i = 1,2,\ldots,4n-1$. By direct calculations, we have $n_i=(i+1)2^i$.
It is both consistent and convenient to define $w_0 = u_0 =n_0 =1$.

\begin{lem} $\sum\nolimits_{i=2}^{4n+1}\frac{1}{\rho_{i}}= n(4n+2)/3$.
\end{lem}

\noindent{\bf Proof.} Let $L_A[i]$ denote the principal submatrix of $L_A[i]$ by deleting its $i$-th row and column. Then
$$(-1)^{4n}a_{4n} = \sum_{i=1}^{4n+1}det L_A[i]= \sum_{i=1}^{4n+1}det \begin{pmatrix} W_{i-1} &0\\ 0 & U_{4n+1-i}\end{pmatrix} =\sum_{i=1}^{4n+1}w_{i-1} w_{4n+1-i}=(4n+1)2^{4n}.$$

On the other hand, let $L_A[i,j]$ denote the submatrix of $L_A[i]$ by deleting its $i$-th and $j$-th rows and columns. Then
\begin{align*}
(-1)^{4n-1}a_{4n-1} &=\sum_{1\leq{i}<j\leq4n+1}det L_A[i,j] \\
   &=\sum_{1\leq{i}<j\leq 4n+1}det\begin{pmatrix} W_{i-1}& & \\& N_{j-i-1} &\\& &U_{4n+1-j}\end{pmatrix}\\
   &=\sum_{1\leq{i}<j\leq 4n+1} n_{j-i-1} w_{i-1} u_{4n+1-j}\\
   &=\sum_{1\leq{i}<j\leq 4n+1} (j-i)2^{j-i-1} 2^{i-1} 2^{4n+1-j}\\
   &=\frac{4n(4n+1)(4n+2)2^{4n-2} }{3},
\end{align*}
\noindent and the result follows from Equation (3.3). \hfill\rule{1ex}{1ex}

Now, from Equations (3.1) -- (3.2) and Lemma 3.1, we immediately have the following result.

\begin{thm} $Kf(Q_{n}) = (32n^{3}+44n^{2}+17n+2) /3$.
\end{thm}

The values of $Kf(Q_{n})$ for $1 \leq n \leq 15$ are listed in Table 1.
\begin{table}[htbp]
\setlength{\abovecaptionskip}{0.05cm} \centering\vspace{.3cm}
\caption{The values of $Kf(Q_{n})$, $1 \leq n \leq 15$.}
\begin{tabular}{c|c|c|c|c|c|c|c|c|c}
  \hline
  $n$ & $Kf(Q_{n})$ & $n$ & $Kf(Q_{n})$ & $n$ & $Kf(Q_{n})$ & $n$ & $Kf(Q_{n})$ & $n$ & $Kf(Q_{n})$ \\
  \hline
  $1$ & $31.67$ & $4$ & $940.67$ & $7$ & $4417.67$ & $10$ & $12190.70$ & $13$ & $25987.70$ \\
  $2$ & $156.00$ & $5$ & $1729.00$ & $8$ & $6446.00$ & $11$ & $16035.00$ & $14$ & $32224.00$\\
  $3$ & $437.67$ & $6$ & $2866.67$ & $9$ & $9015.67$ & $12$ & $20612.70$ & $15$ & $39385.70$\\
  \hline
\end{tabular}
\end{table}

\begin{thm} $\lim\nolimits_{n\rightarrow\infty} Kf(Q_{n})/W(Q_{n})=1/4$.
\end{thm}

\noindent{\bf Proof.} For a subset $U \subseteq V$, let $W(U) = \Sigma_{i \in U}\Sigma_{j \in V} d_{ij}$, i.e., the sum of the sum of the distances between each vertex in $U$ and other vertices in $Q_n$.

(1) Let $U_1=\{1,1',4n+1,(4n+1)'\}$, then $W(U_1) =4(16n^{2}+4n+1)$;

(2) Let $U_2 =\{4i-2,(4i-2)'|i = 1,2,\ldots,n\}$, then
$$W(U_2)=2\sum_{i=1}^{n}(1+\sum_{j=1}^{4i-3}j+\sum_{j=1}^{4n-4i+3}j)=\frac{4}{3}(8n^{3}+3n^{2}+n);$$

(3) Let $U_3 = \{4i-1,(4i-1)'|i = 1,2,\ldots,n\}$, then
$$W(U_3) =2\sum_{i=1}^{n}(1+\sum_{j=1}^{4i-2}j+\sum_{j=1}^{4n-4i+2}j)=\frac{2}{3}(16n^{3}+6n^{2}-n);$$

(4) Let $U_4 = \{4i,(4i)'|i = 1,2,\ldots,n\}$, then
$$W(U_4) =2\sum_{i=1}^{n}(1+\sum_{j=1}^{4i-1}j+\sum_{j=1}^{4n-4i+1}j)=\frac{4}{3}(8n^{3}+3n^{2}+n);$$

(5) Let $U_5 = \{4i+1,(4i+1)'|i = 1,2,\ldots,n-1\}$, then
$$W(U_5) =2\sum_{i=1}^{n-1}(1+\sum_{j=1}^{4i}j+\sum_{j=1}^{4n-4i}j)=\frac{2}{3}(16n^{3}-18n^{2}+5n-3).$$

\noindent Hence $W(Q_{n}) =(W(U_1)+W(U_2)+W(U_3)+W(U_4)+W(U_5))/2=(128n^{3}+96n^{2}+40n)/3$, and from Theorem 3.2 we immediately have
$\lim_{n\rightarrow\infty} Kf(Q_{n})/W(Q_{n})=1/4$. \hfill\rule{1ex}{1ex}

\begin{thm} $\tau(Q_n)=2^{10n+2}\cdot3^{2n-1}$.
\end{thm}

\noindent{\bf Proof.} From Lemma 2.2 and the proof of Lemma 3.1, we have
\begin{align*}
\tau(Q_{n}) &=\frac{1}{8n+2} \cdot \prod_{i=2}^{4n+1}\rho_i \cdot \prod_{i=1}^{4n+1}\sigma_i\\
            &=\frac{1}{8n+2} \cdot (-1)^{4n}a^{4n} \cdot 4^{2n+2} \cdot 6^{2n-1} \\\
            &=2^{10n+2}\cdot3^{2n-1}.
\end{align*}

The proof is completed. \hfill\rule{1ex}{1ex}

The values of $\tau(Q_{n})$, $1 \leq n \leq 8$, are shown in Table 3.
\begin{table}[htbp]
\setlength{\abovecaptionskip}{0.05cm}
 \centering \vspace{.3cm}
\caption{The number of spanning trees of $Q_n$, $1 \leq n \leq 8$.}
\begin{tabular}{c|c|c|c}
  \hline
  $n$ & $\tau(Q_{n})$ & $n$ & $\tau(Q_{n})$ \\
  \hline
  $1$ & $12288$ & $5$ & $88644351465533472768$ \\
  $2$ & $113246208$ & $6$ & $816946343106356485029888$ \\
  $3$ & $1043677052928$ & $7$ & $7528977498068181366035447808$ \\
  $4$ & $9618527719784448$ & $8$ & $69387056622196359469382686998528$ \\
  \hline
\end{tabular}
\end{table}

\section{Degree-Kirchhoff index and the complexity of $Q_n$}

By the same method used in Section 3, we will determine the degree-Kirchhoff index of $Q_n$ in this section.
It is easily seen that
\begin{eqnarray*}
\mathcal{L}_{V_1 V_1}&=&
\left(
  \begin{array}{ccccccccccc}
    1 & \frac{-1}{\sqrt{12}} & & & & & & & & &\\
    \frac{-1}{\sqrt{12}} & 1 & \frac{-1}{4} & & & & & & & &\\
    & \frac{-1}{4} & 1 & \frac{-1}{\sqrt{20}} & & & & & & &\\
    & & \frac{-1}{\sqrt{20}} & 1 & \frac{-1}{5} & & & & & &\\
    & & & \frac{-1}{5} & 1 & \frac{-1}{\sqrt{20}} & & & & &\\
    & & & & \frac{-1}{\sqrt{20}} & 1 & \frac{-1}{4} & & & &\\
    & & & & & & \ddots & & & &\\
    & & & & & & \frac{-1}{\sqrt{20}} & 1 & \frac{-1}{4} & &\\
    & & & & & & & \frac{-1}{4} & 1 & \frac{-1}{\sqrt{20}} &\\
    & & & & & & & & \frac{-1}{\sqrt{20}} & 1 & \frac{-1}{\sqrt{15}}\\
    & & & & & & & & & \frac{-1}{\sqrt{15}} & 1\\
  \end{array}
\right)_{(4n+1)\times(4n+1)}
\end{eqnarray*}

\noindent and

\begin{eqnarray*}
\mathcal{L}_{V_1 V_2}&=&
\left(
  \begin{array}{ccccccccccc}
    \frac{-1}{3} & \frac{-1}{\sqrt{12}} & & & & & & & & &\\
    \frac{-1}{\sqrt{12}} & 0 & \frac{-1}{4} & & & & & & & &\\
    & \frac{-1}{4} & 0 & \frac{-1}{\sqrt{20}} & & & & & & &\\
    & & \frac{-1}{\sqrt{20}} & \frac{-1}{5} & \frac{-1}{5} & & & & & &\\
    & & & \frac{-1}{5} & \frac{-1}{5} & \frac{-1}{\sqrt{20}} & & & & &\\
    & & & & \frac{-1}{\sqrt{20}} & 0 & \frac{-1}{4} & & & &\\
    & & & & & & \ddots & & & &\\
    & & & & & & \frac{-1}{\sqrt{20}} & 0 & \frac{-1}{4} & &\\
    & & & & & & & \frac{-1}{4} & 0 & \frac{-1}{\sqrt{20}} &\\
    & & & & & & & & \frac{-1}{\sqrt{20}} & \frac{-1}{5} & \frac{-1}{\sqrt{15}}\\
    & & & & & & & & & \frac{-1}{\sqrt{15}} & \frac{-1}{3}\\
  \end{array}
\right)_{(4n+1)\times(4n+1)}.
\end{eqnarray*}
Thus
\begin{eqnarray*}
\mathcal L_A&=&
\left(
  \begin{array}{ccccccccccc}
    \frac{2}{3} & \frac{-1}{\sqrt{3}} & & & & & & & & &\\
    \frac{-1}{\sqrt{3}} & 1 & \frac{-1}{2} & & & & & & & &\\
    & \frac{-1}{2} & 1 & \frac{-1}{\sqrt{5}} & & & & & & &\\
    & & \frac{-1}{\sqrt{5}} & \frac{4}{5} & \frac{-2}{5} & & & & & &\\
    & & & \frac{-2}{5} & \frac{4}{5} & \frac{-1}{\sqrt{5}} & & & & &\\
    & & & & \frac{-1}{\sqrt{5}} & 1 & \frac{-1}{2} & & & &\\
    & & & & & & \ddots & & & &\\
    & & & & & & \frac{-1}{\sqrt{5}} & 1 & \frac{-1}{2} & &\\
    & & & & & & & \frac{-1}{2} & 1 & \frac{-1}{\sqrt{5}} &\\
    & & & & & & & & \frac{-1}{\sqrt{5}} & \frac{4}{5} & \frac{-2}{\sqrt{15}}\\
    & & & & & & & & & \frac{-2}{\sqrt{15}} & \frac{2}{3}\\
  \end{array}
\right)_{(4n+1)\times(4n+1)}
\end{eqnarray*}

\noindent and $\mathcal L_S = diag(\frac{4}{3},1,1,\frac{6}{5},\frac{6}{5},1,\ldots,1,1,\frac{6}{5},\frac{4}{3})$.

Suppose $Sp(\mathcal L_A) = \{0=\lambda_{1}, \lambda_{2},\ldots,\lambda_{4n+1}\}$ and
$Sp(\mathcal L_S) = \{\varphi_1, \varphi_2, \ldots, \varphi_{4n+1}\}$.
Then from Lemma 1.2 we have

\begin{equation}
Kf^{*}(Q_{n})=2(18n+1)\Big(\sum_{i=2}^{4n+1}\frac{1}{\lambda_{i}}+\sum_{i=1}^{4n+1}\frac{1}{\varphi_{i}}\Big).
\end{equation}

Obviously, \begin{equation}
\sum_{i=1}^{4n+1}\frac{1}{\varphi_{i}}=1\times{2n}+\frac{5}{6}\times(2n-1)+\frac{3}{4}\times2=\frac{11n+2}{3}.
\end{equation}

To get the value of $\sum_{i=2}^{4n+1}\frac{1}{\lambda_{i}}$,
suppose $\Phi_{\mathcal L_A}(x) =det(xI-\mathcal L_A) =x^{4n+1}+b_{1}x^{3n}+\cdots+b_{4n}x$.
Similarly to Equation (3.3),
we have
\begin{equation}
\sum_{i=2}^{4n+1}\frac{1}{\lambda_{i}}=\frac{(-1)^{4n-1}b_{4n-1}}{(-1)^{4n}b_{4n}}.
\end{equation}

Let $X_{i}$ be the $i$-th order principal submatrix of $\mathcal L_{A}$ formed by its first $i$ rows and $i$ columns, and $x_i = det X_i$, $i=1,2,\ldots,4n$.
By direct calculations we have $x_{1}=\frac{2}{3}$, $x_{2}=\frac{1}{3}$, $x_{3}=\frac{1}{6}$, $x_{4}=\frac{1}{15}$, $x_{5}=\frac{2}{75}$, $x_{6}=\frac{1}{75}$, $x_{4n}=\frac{5}{3} \cdot (\frac{1}{25})^n$,
and for $1\leq k \leq n-1$, it holds that
\begin{eqnarray*}
\begin{cases}
x_{4k}=\frac{4}{5}x_{4k-1}-\frac{1}{5}x_{4k-2}\\
x_{4k+1}=\frac{4}{5}x_{4k}-\frac{4}{25}x_{4k-1}\\
x_{4k+2}=x_{4k+1}-\frac{1}{5}x_{4k}\\
x_{4k+3}=x_{4k+2}-\frac{1}{4}x_{4k+1}
\end{cases}.
\end{eqnarray*}
By solving the above recursive equations we have
\begin{eqnarray*}
\begin{cases}
x_{4k}=\frac{5}{3}\cdot(\frac{1}{25})^{k}\\
x_{4k+1}=\frac{2}{3}\cdot(\frac{1}{25})^{k}\\
x_{4k+2}=\frac{1}{3}\cdot(\frac{1}{25})^{k}\\
x_{4k+3}=\frac{1}{6}\cdot(\frac{1}{25})^{k}\\
\end{cases},~1\leq k \leq n-1.
\end{eqnarray*}

Similarly, let $Y_{i}$ be the $i$-th order principal submatrix of $\mathcal L_{A}$ formed by its last $i$ rows and $i$ columns, and $y_i = det Y_i$, $i=1,2,\ldots,4n$. Then $y_{1}=\frac{2}{3}$, $y_{2}=\frac{4}{15}$, $y_{3}=\frac{2}{15}$, $y_{4}=\frac{1}{15}$, $y_{5}=\frac{2}{75}$, $y_{6}=\frac{4}{375}$, $y_{4n}=\frac{5}{3}\cdot(\frac{1}{25})^{n}$, and for $1\leq k \leq n-1$, it holds that
\begin{eqnarray*}
\begin{cases}
y_{4k}=y_{4k-1}-\frac{1}{4}y_{4k-2}\\
y_{4k+1}=\frac{4}{5}y_{4k}-\frac{1}{5}y_{4k-1}\\
y_{4k+2}=\frac{4}{5}y_{4k+1}-\frac{4}{25}y_{4k}\\
y_{4k+3}=y_{4k+2}-\frac{1}{5}y_{4k+1}
\end{cases},
\end{eqnarray*}
and so
\begin{eqnarray*}
\begin{cases}
y_{4k}=\frac{5}{3}\cdot(\frac{1}{25})^{k}\\
y_{4k+1}=\frac{2}{3}\cdot(\frac{1}{25})^{k}\\
y_{4k+2}=\frac{4}{15}\cdot(\frac{1}{25})^{k}\\
y_{4k+3}=\frac{2}{15}\cdot(\frac{1}{25})^{k}
\end{cases}.
\end{eqnarray*}

Define $x_0=y_0=1$.

\begin{lem} $\sum_{i=2}^{4n+1}\frac{1}{\lambda_i}= \frac{2n(54n^{2}+9n+4)}{18n+1}$. 
\end{lem}

\noindent{\bf Proof.} Let $\mathcal L_A[i]$ denote the principal submatrix of $\mathcal L_A[i]$ by deleting its $i$-th row and column. Then
\begin{align*}
(-1)^{4n}b_{4n}&=\sum_{j=1}^{4n+1}det\mathcal L_A[j] =x_{4n}+y_{4n}+\sum_{j=2}^{4n}det\mathcal L_A[j]\\
&=x_{4n}+y_{4n}+\sum_{k=1}^{n}x_{4(k-1)+3}y_{4(n-k)+1}+\sum_{k=1}^{n-1}x_{4k}y_{4(n-k)}\\
&~~+\sum_{k=0}^{n-1}x_{4k+1}y_{4(n-k-1)+3}+\sum_{k=0}^{n-1}x_{4k+2}y_{4(n-k-1)+2}\\
&=\frac{10}{3}\cdot(\frac{1}{25})^{n}+\sum_{k=1}^{n}\frac{1}{6}\cdot(\frac{1}{25})^{k-1}\cdot\frac{2}{3}\cdot(\frac{1}{25})^{n-k}
+\sum_{k=1}^{n-1}\frac{5}{3}\cdot(\frac{1}{25})^{k}\cdot\frac{5}{3}\cdot(\frac{1}{25})^{n-k}\\
&~~+\sum_{k=0}^{n-1}\frac{2}{3}\cdot\left(\frac{1}{25}\right)^{k}\cdot\frac{2}{15}\cdot(\frac{1}{25})^{n-k-1}+\sum_{k=0}^{n-1}\frac{1}{3}\cdot(\frac{1}{25})^{k}\cdot\frac{4}{15}\cdot(\frac{1}{25})^{n-k-1}\\
&=\frac{18n+1}{45}(\frac{1}{25})^{n-1}.
\end{align*}

From Equation (4.3), it suffices to show
$(-1)^{4n-1}b_{4n-1}=\frac{2(54n^{3}+9n^{2}+4n)}{45}(\frac{1}{25})^{n-1}
$.

Let $\mathcal L_A[i,j]$ denote the submatrix of $\mathcal L_A$ by deleting its $i$-th and $j$-th rows and columns, $1 \leq i < j \leq 4n+1$. Then
\begin{align}\label{sum}
(-1)^{4n-1}b_{4n-1} &=\sum_{1\leq i < j \leq 4n+1}det \mathcal L_A[i,j]
   =\sum_{1\leq i<j \leq 4n+1}det\begin{pmatrix} X_{i-1}& & \\& Z_{i,j} &\\& &Y_{4n+1-j}\end{pmatrix}\\ \nonumber
   &= \sum_{1\leq i < j \leq 4n+1}x_{i-1} y_{4n+1-j} z_{i,j},
\end{align}
\noindent where $Z_{i,j}$ is the $(j-i-1)$-th order principal submatrix of $\mathcal L_A$ formed by its $i+1, i+2, \ldots, j-1$ rows and columns. Let $z_{i,j} = detZ_{i,j}$,
and defined $z_{i,j}=1$ if $j=i+1$.
It is an easy but tedious task to calculate the values of $z_{i,j}$ for different $i$'s and $j$'s.
In fact, in the Appendix, we calculate the values of $z_{i,j}$ in 16 cases, and have the following result.

\begin{equation}\label{z_ij}
z_{i,j}=
  \begin{cases}
\frac{2}{5}(b-a)(\frac{1}{25})^{b-a-1},&i=4a, j=4b, 1\leq a < b \leq n;\\
(4b-4a+1)(\frac{1}{25})^{b-a}, &i=4a, j=4b+1, 1\leq{a}\leq b\leq{n};\\
\frac{2}{5}(4b-4a+2)(\frac{1}{25})^{b-a}, &i=4a, j=4b+2, 1\leq{a}\leq b\leq{n-1};\\
\frac{1}{5}(4b-4a+3)(\frac{1}{25})^{b-a}, &i=4a, j=4b+3, 1\leq{a}\leq b\leq{n-1};\\
\frac{1}{4}(4b-4a-1)(\frac{1}{25})^{b-a-1}, &i=4a+1, j=4b, 0\leq{a}< b\leq{n};\\
\frac{2}{5}(b-a)(\frac{1}{25})^{b-a-1}, &i=4a+1, j=4b+1, 0\leq{a}< b\leq{n};\\
(4b-4a+1)(\frac{1}{25})^{b-a}, &i=4a+1, j=4b+2, 0\leq{a}\leq b\leq{n-1};\\
(2b-2a+1)(\frac{1}{25})^{b-a}, &i=4a+1, j=4b+3, 0\leq{a}\leq b\leq{n-1};\\
(2b-2a-1)(\frac{1}{25})^{b-a-1}, &i=4a+2, j=4b, 0\leq{a}< b\leq{n};\\
\frac{1}{5}(4b-4a-1)(\frac{1}{25})^{b-a-1}, &i=4a+2, j=4b+1, 0\leq{a}< b\leq{n};\\
\frac{2}{25}(4b-4a)(\frac{1}{25})^{b-a-1}, &i=4a+2,j=4b+2, 0\leq{a}< b\leq{n-1};\\
(4b-4a+1)(\frac{1}{25})^{b-a}, &i=4a+2, j=4b+3, 0\leq{a}\leq b\leq{n-1};\\
(4b-4a-3)(\frac{1}{25})^{b-a-1}, &i=4a+3, j=4b, 0\leq{a}<b\leq{n};\\
\frac{2}{5}(4b-4a-2)(\frac{1}{25})^{b-a-1}, &i=4a+3, j=4b+1, 0\leq{a}<b\leq{n};\\
\frac{4}{25}(4b-4a-1)(\frac{1}{25})^{b-a-1}, &i=4a+3, j=4b+2, 0\leq{a}<b\leq{n-1};\\
\frac{2}{25}(4b-4a)(\frac{1}{25})^{b-a-1}, &i=4a+3, j=4b+3, 0\leq{a}<b\leq{n-1}.
\end{cases}
\end{equation}

Accordingly, in the Appendix, we obtain the values of $K_{i,j}$, $0\leq i,j\leq 3$, as follows.
\begin{equation}\label{K_ij}
  \begin{cases}
K_{0,0}=\sum_{1\leq{a}<b\leq{n}} det\mathcal L_A[4a,4b]=
  \frac{5(n^3-n)}{27}\left(\frac{1}{25}\right)^{n-1};\\
K_{0,1}=\sum_{1\leq{a}\leq b\leq n} det\mathcal L_A[4a,4b+1]
  =\frac{20n^3-9n^2+7n}{108}\left(\frac{1}{25}\right)^{n-1};\\
K_{0,2}=\sum_{1\leq{a}\leq b\leq{n-1}} det\mathcal L_A[4a,4b+2]
  =\frac{2(2n^3-3n^{2}+n)}{27}\left(\frac{1}{25}\right)^{n-1};\\
K_{0,3}=\sum_{1\leq{a}\leq b\leq{n-1}} det\mathcal L_A[4a,4b+3]
  =\frac{4n^3-3n^{2}-n}{27}\left(\frac{1}{25}\right)^{n-1};\\
K_{1,0}=\sum_{0\leq{a}<b\leq{n}} det\mathcal L_A[4a+1,4b]
  =\frac{20n^3+21n^2+13n}{108} \left(\frac{1}{25}\right)^{n-1};\\
K_{1,1}=\sum_{0\leq{a}<b\leq{n}} det\mathcal L_A[4a+1,4b+1]
  =\frac{25n^3+15n^2+14n}{135}\left(\frac{1}{25}\right)^{n-1};\\
K_{1,2}=\sum_{0\leq{a}\leq b\leq{n-1}} det\mathcal L_A[4a+1,4b+2]
  =\frac{20n^3-9n^2+7n}{135}\left(\frac{1}{25}\right)^{n-1};\\
K_{1,3}=\sum_{0\leq{a}\leq b\leq{n-1}} det\mathcal L_A[4a+1,4b+3]=
  \frac{20n^3+6n^2+10n}{135}\left(\frac{1}{25}\right)^{n-1};\\
K_{2,0}=\sum_{0\leq{a}<b\leq{n}} det\mathcal L_A[4a+2,4b]=
  \frac{2(2n^3+3n^{2}+n)}{27}\left(\frac{1}{25}\right)^{n-1};\\
K_{2,1}=\sum_{0\leq{a}<b\leq n} det\mathcal L_A[4a+2,4b+1]=
  \frac{20n^3+21n^2+13n}{135}\left(\frac{1}{25}\right)^{n-1};\\
K_{2,2}=\sum_{0\leq{a}<b\leq{n-1}} det\mathcal L_A[4a+2,4b+2]=
  \frac{16(n^3-n)}{135}\left(\frac{1}{25}\right)^{n-1};\\
K_{2,3}=\sum_{0\leq{a}<b\leq{n-1}} det\mathcal L_A[4a+2,4b+3]=
  \frac{4(4n^3+3n^{2}-n)}{135}\left(\frac{1}{25}\right)^{n-1};\\
K_{3,0}=\sum_{0\leq{a}<b\leq{n}} det\mathcal L_A[4a+3,4b]=
  \frac{4n^3+3n^2-n}{27}\left(\frac{1}{25}\right)^{n-1};\\
K_{3,1}=\sum_{0\leq{a}<b\leq{n}} det\mathcal L_A[4a+3,4b+1]=
  \frac{20n^3+6n^2+10n}{135}\left(\frac{1}{25}\right)^{n-1};\\
K_{3,2}=\sum_{0\leq{a}<b\leq{n-1}} det\mathcal L_A[4a+3,4b+2]=
   \frac{4(4n^3-3n^{2}-n)}{135}\left(\frac{1}{25}\right)^{n-1};\\
K_{3,3}=\sum_{0\leq{a}<b\leq{n-1}} det\mathcal L_A[4a+3,4b+3]=
  \frac{16(n^3-n)}{135}\left(\frac{1}{25}\right)^{n-1}.
\end{cases}
\end{equation}
Therefore, $
(-1)^{4n-1}b_{4n-1}=\sum_{0\leq i,j\leq 3}K_{i,j} =\frac{2n(54n^{2}+9n+4)}{45}(\frac{1}{25})^{n-1}. \hfill\rule{1ex}{1ex}$

From Equations (4.1)--(4.2) and Lemma 4.1, we immediately have the following result.
\begin{thm} $Kf^{*}(Q_{n})= 2(324n^{3}+252n^{2}+71n+2)/3$.
\end{thm}

\begin{thm} $\lim\nolimits_{n\rightarrow\infty} Kf^{*}(Q_{n})/Gut(Q_{n})=1/4$.
\end{thm}

\noindent{\bf Proof.} For a subset $U \subseteq V$,
let $Gut(U) = \Sigma_{i \in U}\Sigma_{j \in V} d_i d_j d_{ij}$.

(1) Let $U_1=\{1,1'\}$, then
\begin{align*}
Gut(U_1)&=2\left[3\cdot3\cdot(1+4n+4n)+\sum_{i=1}^{n}2\cdot3\cdot5\cdot(4i-1)
  +\sum_{i=1}^{n-1}2\cdot3\cdot5\cdot(4i)\right.\\
&~~~\left.+\sum_{i=0}^{n-1}2\cdot3\cdot4\cdot(4i+1)+\sum_{i=0}^{n-1}2\cdot3\cdot4\cdot(4i+2)\right]\\
&=18(24n^{2}+2n+1).
\end{align*}

(2) Let $U_2=\{4n+1,(4n+1)'\}$, then
\begin{align*}
Gut(U_2)&=2\left[3\cdot3\cdot(1+4n+4n)+\sum_{i=1}^{n}2\cdot3\cdot5\cdot(4n-4i+1)+\sum_{i=1}^{n-1}2\cdot3\cdot4\cdot(4n-4i+2)\right.\\
&~~~\left.+\sum_{i=1}^{n}2\cdot3\cdot4\cdot(4n-4i+3)+\sum_{i=1}^{n}2\cdot3\cdot5\cdot(4n-4i+4)\right]\\
&=6(72n^{2}+42n-13).
\end{align*}

(3) Let $U_3 =\{4i-2,(4i-2)'|i = 1,2,\ldots,n\}$, then
\begin{align*}
Gut(U_3)&=2\sum_{i=1}^{n}\left[4\cdot4\cdot2+2\cdot3\cdot4\cdot(4i-3)+2\cdot3\cdot4\cdot(4n-4i+3)+\sum_{k=2}^{i}2\cdot4\cdot5\cdot(4i-4k+1)\right.\\
&~~~+\sum_{k=i+1}^{n}2\cdot4\cdot5\cdot(4k-4i-1)+\sum_{k=1}^{i-1}2\cdot4\cdot4\cdot(4i-4k-1)+\sum_{k=i}^{n}2\cdot4\cdot4\cdot(4k-4i+1)\\
&~~~+\sum_{k=1}^{i-1}2\cdot4\cdot4\cdot(4i-4k)+\sum_{k=i+1}^{n}2\cdot4\cdot4\cdot(4k-4i)+\sum_{k=1}^{i-1}2\cdot4\cdot5\cdot(4i-4k-2)\\
&~~~\left.+\sum_{k=i}^{n}2\cdot4\cdot5\cdot(4k-4i+2)\right]\\
&=32(12n^{3}+n^{2}+2n).
\end{align*}

(4) If $U_4 = \{4i-1,(4i-1)'|i = 1,2,\ldots,n\}$, then
\begin{align*}
Gut(U_4)&=2\sum_{i=1}^{n}\left[2\cdot4\cdot4\cdot2
+2\cdot3\cdot4\cdot(4i-2)+2\cdot3\cdot4\cdot(4n-4i+2)+\sum_{k=1}^{i-1}2\cdot4\cdot4\cdot(4i-4k)\right.\\
&~~~+\sum_{k=i+1}^{n}2\cdot4\cdot4\cdot(4k-4i)+\sum_{k=1}^{i}2\cdot4\cdot4\cdot(4i-4k+1)+\sum_{k=i+1}^{n}2\cdot4\cdot4\cdot(4k-4i-1)\\
&~~~+\sum_{k=2}^{i}2\cdot4\cdot5\cdot(4i-4k+2)+\sum_{k=i+1}^{n}2\cdot4\cdot5\cdot(4k-4i-2)+\sum_{k=1}^{i-1}2\cdot4\cdot5\cdot(4i-4k-1)\\
&~~~\left.+\sum_{k=i}^{n}2\cdot4\cdot5\cdot(4k-4i+1)\right]\\
&=16(24n^{3}+2n^{2}-n).
\end{align*}

(5) If $U_5 = \{4i,(4i)'|i = 1,2,\ldots,n\}$, then
\begin{align*}
Gut(U_5)&=2\sum_{i=1}^{n}\left[5\cdot5\cdot1+2\cdot3\cdot5\cdot(4i-1)
+2\cdot3\cdot5\cdot(4n-4i+1)+\sum_{k=1}^{i}2\cdot4\cdot5\cdot(4i-4k+1)\right.\\
&~~~+\sum_{k=i+1}^{n}2\cdot4\cdot5\cdot(4k-4i-1)+\sum_{k=1}^{i}2\cdot4\cdot5\cdot(4i-4k+2)+\sum_{k=i+1}^{n}2\cdot4\cdot5\cdot(4k-4i-2)\\
&~~~+\sum_{k=2}^{i}2\cdot5\cdot5\cdot(4i-4k+3)+\sum_{k=i+1}^{n}2\cdot5\cdot5\cdot(4k-4i-3)+\sum_{k=1}^{i-1}2\cdot5\cdot5\cdot(4i-4k)\\
&~~~\left.+\sum_{k=i+1}^{n}2\cdot5\cdot5\cdot(4k-4i)\right]\\
&=10(48n^{3}+4n^{2}+n).
\end{align*}

(6) If $U_6 = \{4i-3,(4i-3)'|i = 2,\ldots,n\}$, then
\begin{align*}
Gut(U_6)&=2\sum_{i=2}^{n}\left[5\cdot5\cdot1+2\cdot3\cdot5\cdot(4i-4)+2\cdot3\cdot5\cdot(4n-4i+4)+\sum_{k=1}^{i-1}2\cdot4\cdot5\cdot(4i-4k-2)\right.\\
&~~~+\sum_{k=i}^{n}2\cdot4\cdot5\cdot(4k-4i+2)+\sum_{k=2}^{i-1}2\cdot4\cdot5\cdot(4i-4k-1)+\sum_{k=i}^{n}2\cdot4\cdot5\cdot(4k-4i+1)\\
&~~~+\sum_{k=2}^{i-1}2\cdot5\cdot5\cdot(4i-4k)+\sum_{k=i+1}^{n}2\cdot5\cdot5\cdot(4k-4i)+\sum_{k=1}^{i-1}2\cdot5\cdot5\cdot(4i-4k-3)\\
&~~~\left.+\sum_{k=i}^{n}2\cdot5\cdot5\cdot(4k-4i+3)\right]\\
&=10(48n^{3}-84n^{2}+49n-13).
\end{align*}

Hence $Gut(Q_{n}) =(Gut(U_1)+Gut(U_2)+ \cdots +Gut(U_6))/2=864n^{3}+64n^{2}+418n-95$, and from Theorem 4.2 we immediately have
$\lim_{n\rightarrow\infty} Kf^{*}(Q_{n})/Gut(Q_{n})=1/4$. \hfill\rule{1ex}{1ex}

The values of $Kf^{*}(Q_{n})$ for $1 \leq n \leq 15$ are listed in Table 2.
\begin{table}[htbp]
\setlength{\abovecaptionskip}{0.05cm} \centering\vspace{.3cm}
\caption{The values of $Kf^*(Q_n)$, $1 \leq n \leq 15$.}
\begin{tabular}{c|c|c|c|c|c|c|c|c|c}
  \hline
  $n$ & $Kf^{*}(Q_{n})$ & $Q_{n}$ & $Kf^{*}(Q_{n})$ & $n$ & $Kf^{*}(Q_{n})$ & $n$ & $Kf^{*}(Q_{n})$ & $Q_{n}$ & $Kf^{*}(Q_{n}))$ \\
  \hline
  $1$ & $432.67$ & $4$ & $16702.70$ & $7$ & $82652.70$ & $10$ & $233274.67$ & $13$ & $503560.67$ \\
  $2$ & $2496.00$ & $5$ & $31438.00$ & $8$ & $121724.00$ & $11$ & $308316.00$ & $14$ & $626296.00$
  \\
  $3$ & $7487.33$ & $6$ & $52989.30$ & $9$ & $171499.34$ & $12$ & $398009.34$ & $15$ & $767511.34$
  \\
  \hline
\end{tabular}
\end{table}

%

\section*{Appendix: The values of $z_{j-i-1}$'s for different $i$'s and $j$'s}

For $1 \leq i < j \leq 4n+1$, the values of $z_{i,j} = det Z_{i,j}$ (listed in Eq. (\ref{z_ij})) can be easily obtained.
Since it is defined that $z_{i,j} = 1$ if $j=i+1$, here we assume $j>i+1$.
We distinguish the following 16 cases.

{\bf Case 1.} $i=4a$ and $j=4b$, $1\leq a < b \leq n$. Then
\begin{align*}
z_{i,j}&=
\left|
  \begin{array}{cccccccc}
    \frac{4}{5} & -\frac{1}{\sqrt{5}} & & & & & &\\
    -\frac{1}{\sqrt{5}} & 1 & -\frac{1}{2} & & & & &\\
    & -\frac{1}{2} & 1 & -\frac{1}{\sqrt{5}} & & & &\\
    & & & \ddots & & & &\\
    & & &  & \frac{4}{5} & -\frac{2}{5} & &\\
    & & & & -\frac{2}{5} & \frac{4}{5} & -\frac{1}{\sqrt{5}} &\\
    & & & & & -\frac{1}{\sqrt{5}} & 1 & -\frac{1}{2}\\
    & & & & & & -\frac{1}{2} & 1\\
  \end{array}
\right|_{(4b-4a-1)\times (4b-4a-1)}\\
&=\frac{2}{5}(b-a)\left(\frac{1}{25}\right)^{b-a-1}.
\end{align*}

{\bf Case 2.} $i=4a$ and $j=4b+1$, $1\leq{a}\leq b\leq n$. Then
\begin{align*}
z_{i,j}&= \left|
  \begin{array}{cccccccc}
    \frac{4}{5} & -\frac{1}{\sqrt{5}} & & & & & &\\
    -\frac{1}{\sqrt{5}} & 1 & -\frac{1}{2} & & & & &\\
    & -\frac{1}{2} & 1 & -\frac{1}{\sqrt{5}} & & & &\\
    & & & \ddots & & & &\\
    & & &  & \frac{4}{5} & -\frac{1}{\sqrt{5}} & &\\
    & & & & -\frac{1}{\sqrt{5}} & 1 & -\frac{1}{2} &\\
    & & & & & -\frac{1}{2} & 1 &  -\frac{1}{\sqrt{5}}\\
    & & & & & & -\frac{1}{\sqrt{5}} & \frac{4}{5}\\
  \end{array}
\right|_{(4b-4a)\times (4b-4a)}\\
&=(4b-4a+1)\left(\frac{1}{25}\right)^{b-a}.
\end{align*}

{\bf Case 3.} $i=4a$ and $j=4b+2$, $1\leq{a}\leq b\leq{n-1}$. Then
\begin{align*}
z_{i,j}&= \left|
  \begin{array}{cccccccc}
    \frac{4}{5} & -\frac{1}{\sqrt{5}} & & & & & &\\
    -\frac{1}{\sqrt{5}} & 1 & -\frac{1}{2} & & & & &\\
    & -\frac{1}{2} & 1 & -\frac{1}{\sqrt{5}} & & & &\\
    & & & \ddots & & & &\\
    & & &  & 1 & -\frac{1}{2} & &\\
    & & & & -\frac{1}{2} & 1 & -\frac{1}{\sqrt{5}} &\\
    & & & & & -\frac{1}{\sqrt{5}} & \frac{4}{5} &  -\frac{2}{5}\\
    & & & & & & -\frac{2}{5} & \frac{4}{5}\\
  \end{array}
\right|_{(4b-4a+1)\times (4b-4a+1)}\\
&=\frac{2}{5}(4b-4a+2)\left(\frac{1}{25}\right)^{b-a}.
\end{align*}

{\bf Case 4.} $i=4a$ and $j=4b+3$, $1\leq{a}\leq b\leq{n-1}$. Then
\begin{align*}
z_{i,j}&= \left|
  \begin{array}{cccccccc}
    \frac{4}{5} & -\frac{1}{\sqrt{5}} & & & & & &\\
    -\frac{1}{\sqrt{5}} & 1 & -\frac{1}{2} & & & & &\\
    & -\frac{1}{2} & 1 & -\frac{1}{\sqrt{5}} & & & &\\
    & & & \ddots & & & &\\
    & & &  & 1 & -\frac{1}{\sqrt{5}} & &\\
    & & & & -\frac{1}{\sqrt{5}} & \frac{4}{5} & -\frac{2}{5} &\\
    & & & & & -\frac{2}{5} & \frac{4}{5} &  -\frac{1}{\sqrt{5}}\\
    & & & & & & -\frac{1}{\sqrt{5}} & 1\\
  \end{array}
\right|_{(4b-4a+2)\times (4b-4a+2)}\\
&=\frac{1}{5}(4b-4a+3)\left(\frac{1}{25}\right)^{b-a}.
\end{align*}

{\bf Case 5.} $i=4a+1$ and $j=4b$, $0\leq{a}< b\leq{n}$. Then
\begin{align*}
z_{i,j}&= \left|
  \begin{array}{cccccccc}
    1 & -\frac{1}{2} & & & & & &\\
    -\frac{1}{2} & 1 & -\frac{1}{\sqrt{5}} & & & & &\\
    & -\frac{1}{\sqrt{5}} & \frac{4}{5} & -\frac{2}{5} & & & &\\
    & & & \ddots & & & &\\
    & & &  & \frac{4}{5} & -\frac{2}{5} & &\\
    & & & & -\frac{2}{5} & \frac{4}{5} & -\frac{1}{\sqrt{5}} &\\
    & & & & & -\frac{1}{\sqrt{5}} & 1 &  -\frac{1}{2}\\
    & & & & & & -\frac{1}{2} & 1\\
  \end{array}
\right|_{(4b-4a-2)\times (4b-4a-2)}\\
&=\frac{1}{4}(4b-4a-1)\left(\frac{1}{25}\right)^{b-a-1}.
\end{align*}

{\bf Case 6.} $i=4a+1$ and $j=4b+1$, $0\leq{a}< b\leq n$. Then
\begin{align*}
z_{i,j}&= \left|
  \begin{array}{cccccccc}
    1 & -\frac{1}{2} & & & & & &\\
    -\frac{1}{2} & 1 & -\frac{1}{\sqrt{5}} & & & & &\\
    & -\frac{1}{\sqrt{5}} & \frac{4}{5} & -\frac{2}{5} & & & &\\
    & & & \ddots & & & &\\
    & & &  & \frac{4}{5} & -\frac{1}{\sqrt{5}} & &\\
    & & & & -\frac{1}{\sqrt{5}} & 1 & -\frac{1}{2} &\\
    & & & & & -\frac{1}{2} & 1 &  -\frac{1}{\sqrt{5}}\\
    & & & & & & -\frac{1}{\sqrt{5}} & \frac{4}{5}\\
  \end{array}
\right|_{(4b-4a-1)\times (4b-4a-1)}\\
&=\frac{2}{5}(b-a)\left(\frac{1}{25}\right)^{b-a-1}.
\end{align*}

{\bf Case 7.} $i=4a+1$ and $j=4b+2$, $0\leq{a}\leq b\leq{n-1}$. Then
\begin{align*}
z_{i,j}&= \left|
  \begin{array}{cccccccc}
    1 & -\frac{1}{2} & & & & & &\\
    -\frac{1}{2} & 1 & -\frac{1}{\sqrt{5}} & & & & &\\
    & -\frac{1}{\sqrt{5}} & \frac{4}{5} & -\frac{2}{5} & & & &\\
    & & & \ddots & & & &\\
    & & &  & 1 & -\frac{1}{2} & &\\
    & & & & -\frac{1}{2} & 1 & -\frac{1}{\sqrt{5}} &\\
    & & & & & -\frac{1}{\sqrt{5}} & \frac{4}{5} &  -\frac{2}{5}\\
    & & & & & & -\frac{2}{5} & \frac{4}{5}\\
  \end{array}
\right|_{(4b-4a)\times (4b-4a)}\\
&=(4b-4a+1)\left(\frac{1}{25}\right)^{b-a}.
\end{align*}

{\bf Case 8.} $i=4a+1$ and $j=4b+3$, $0\leq{a}\leq b\leq{n-1}$. Then
\begin{align*}
z_{i,j}&= \left|
  \begin{array}{cccccccc}
    1 & -\frac{1}{2} & & & & & &\\
    -\frac{1}{2} & 1 & -\frac{1}{\sqrt{5}} & & & & &\\
    & -\frac{1}{\sqrt{5}} & \frac{4}{5} & -\frac{2}{5} & & & &\\
    & & & \ddots & & & &\\
    & & &  & 1 & -\frac{1}{\sqrt{5}} & &\\
    & & & & -\frac{1}{\sqrt{5}} & \frac{4}{5} & -\frac{2}{5} &\\
    & & & & & -\frac{2}{5} & \frac{4}{5} &  -\frac{1}{\sqrt{5}}\\
    & & & & & & -\frac{1}{\sqrt{5}} & 1\\
  \end{array}
\right|_{(4b-4a+1)\times (4b-4a+1)}\\
&=(2b-2a+1)\left(\frac{1}{25}\right)^{b-a}.
\end{align*}

{\bf Case 9.} $i=4a+2$ and $j=4b$, $0\leq{a}< b\leq{n}$. Then
\begin{align*}
z_{i,j}&= \left|
  \begin{array}{cccccccc}
    1 & -\frac{1}{\sqrt{5}} & & & & & &\\
    -\frac{1}{\sqrt{5}} & \frac{4}{5} & -\frac{2}{5} & & & & &\\
    & -\frac{2}{5} & \frac{4}{5} & -\frac{1}{\sqrt{5}} & & & &\\
    & & & \ddots & & & &\\
    & & &  & \frac{4}{5} & -\frac{2}{5} & &\\
    & & & & -\frac{2}{5} & \frac{4}{5} & -\frac{1}{\sqrt{5}} &\\
    & & & & & -\frac{1}{\sqrt{5}} & 1 &  -\frac{1}{2}\\
    & & & & & & -\frac{1}{2} & 1\\
  \end{array}
\right|_{(4b-4a-3)\times (4b-4a-3)}\\
&=(2b-2a-1)\left(\frac{1}{25}\right)^{b-a-1}.
\end{align*}

{\bf Case 10.} $i=4a+2$ and $j=4b+1$, $0\leq{a}< b\leq{n}$. Then
\begin{align*}
z_{i,j}&=
\left|
  \begin{array}{cccccccc}
    1 & -\frac{1}{\sqrt{5}} & & & & & &\\
    -\frac{1}{\sqrt{5}} & \frac{4}{5} & -\frac{2}{5} & & & & &\\
    & -\frac{2}{5} & \frac{4}{5} & -\frac{1}{\sqrt{5}} & & & &\\
    & & & \ddots & & & &\\
    & & &  & \frac{4}{5} & -\frac{1}{\sqrt{5}} & &\\
    & & & & -\frac{1}{\sqrt{5}} & 1 & -\frac{1}{2} &\\
    & & & & & -\frac{1}{2} & 1 &  -\frac{1}{\sqrt{5}}\\
    & & & & & & -\frac{1}{\sqrt{5}} & \frac{4}{5}\\
  \end{array}
\right|_{(4b-4a-2)\times (4b-4a-2)}\\
&=\frac{1}{5}(4b-4a-1)\left(\frac{1}{25}\right)^{b-a-1}.
\end{align*}

{\bf Case 11.} $i=4a+2$ and $j=4b+2$, $0\leq{a}< b\leq{n-1}$. Then
\begin{align*}
z_{i,j}&= \left|
  \begin{array}{cccccccc}
    1 & -\frac{1}{\sqrt{5}} & & & & & &\\
    -\frac{1}{\sqrt{5}} & \frac{4}{5} & -\frac{2}{5} & & & & &\\
    & -\frac{2}{5} & \frac{4}{5} & -\frac{1}{\sqrt{5}} & & & &\\
    & & & \ddots & & & &\\
    & & &  & 1 & -\frac{1}{2} & &\\
    & & & & -\frac{1}{2} & 1 & -\frac{1}{\sqrt{5}} &\\
    & & & & & -\frac{1}{\sqrt{5}} & \frac{4}{5} &  -\frac{2}{5}\\
    & & & & & & -\frac{2}{5} & \frac{4}{5}\\
  \end{array}
\right|_{(4b-4a-1)\times (4b-4a-1)}\\
&=\frac{2}{25}(4b-4a)\left(\frac{1}{25}\right)^{b-a-1}.
\end{align*}

{\bf Case 12.} $i=4a+2$ and $j=4b+3, 0\leq{a}\leq b\leq{n-1}.$ Then
\begin{align*}
z_{i,j}&= \left|
  \begin{array}{cccccccc}
    1 & -\frac{1}{\sqrt{5}} & & & & & &\\
    -\frac{1}{\sqrt{5}} & \frac{4}{5} & -\frac{2}{5} & & & & &\\
    & -\frac{2}{5} & \frac{4}{5} & -\frac{1}{\sqrt{5}} & & & &\\
    & & & \ddots & & & &\\
    & & &  & 1 & -\frac{1}{\sqrt{5}} & &\\
    & & & & -\frac{1}{\sqrt{5}} & \frac{4}{5} & -\frac{2}{5} &\\
    & & & & & -\frac{2}{5} & \frac{4}{5} &  -\frac{1}{\sqrt{5}}\\
    & & & & & & -\frac{1}{\sqrt{5}} & 1\\
  \end{array}
\right|_{(4b-4a)\times (4b-4a)}\\
&=(4b-4a+1)\left(\frac{1}{25}\right)^{b-a}.
\end{align*}

{\bf Case 13.} $i=4a+3$ and $j=4b$, $0\leq{a}<b\leq{n}$. Then
\begin{align*}
z_{i,j}&= \left|
  \begin{array}{cccccccc}
    \frac{4}{5} & -\frac{2}{5} & & & & & &\\
    -\frac{2}{5} & \frac{4}{5} & -\frac{1}{\sqrt{5}} & & & & &\\
    & -\frac{1}{\sqrt{5}} & 1 & -\frac{1}{2} & & & &\\
    & & & \ddots & & & &\\
    & & &  & \frac{4}{5} & -\frac{2}{5} & &\\
    & & & & -\frac{2}{5} & \frac{4}{5} & -\frac{1}{\sqrt{5}} &\\
    & & & & & -\frac{1}{\sqrt{5}} & 1 &  -\frac{1}{2}\\
    & & & & & & -\frac{1}{2} & 1\\
  \end{array}
\right|_{(4b-4a-4)\times (4b-4a-4)}\\
&=(4b-4a-3)\left(\frac{1}{25}\right)^{b-a-1}.
\end{align*}

{\bf Case 14.} $i=4a+3$ and $j=4b+1$, $0\leq{a}<b\leq{n}$. Then
\begin{align*}
z_{i,j}&= \left|
  \begin{array}{cccccccc}
    \frac{4}{5} & -\frac{2}{5} & & & & & &\\
    -\frac{2}{5} & \frac{4}{5} & -\frac{1}{\sqrt{5}} & & & & &\\
    & -\frac{1}{\sqrt{5}} & 1 & -\frac{1}{2} & & & &\\
    & & & \ddots & & & &\\
    & & &  & \frac{4}{5} & -\frac{1}{\sqrt{5}} & &\\
    & & & & -\frac{1}{\sqrt{5}} & 1 & -\frac{1}{2} &\\
    & & & & & -\frac{1}{2} & 1 &  -\frac{1}{\sqrt{5}}\\
    & & & & & & -\frac{1}{\sqrt{5}} & \frac{4}{5}\\
  \end{array}
\right|_{(4b-4a-3)\times (4b-4a-3)}\\
&=\frac{2}{5}(4b-4a-2)\left(\frac{1}{25}\right)^{b-a-1}.
\end{align*}

{\bf Case 15.} $i=4a+3$ and $j=4b+2$, $0\leq{a}<b\leq{n-1}$. Then
\begin{align*}
z_{i,j}&= \left|
  \begin{array}{cccccccc}
    \frac{4}{5} & -\frac{2}{5} & & & & & &\\
    -\frac{2}{5} & \frac{4}{5} & -\frac{1}{\sqrt{5}} & & & & &\\
    & -\frac{1}{\sqrt{5}} & 1 & -\frac{1}{2} & & & &\\
    & & & \ddots & & & &\\
    & & &  & 1 & -\frac{1}{2} & &\\
    & & & & -\frac{1}{2} & 1 & -\frac{1}{\sqrt{5}} &\\
    & & & & & -\frac{1}{\sqrt{5}} & \frac{4}{5} &  -\frac{2}{5}\\
    & & & & & & -\frac{2}{5} & \frac{4}{5}\\
  \end{array}
\right|_{(4b-4a-2)\times (4b-4a-2)}\\
&=\frac{4}{25}(4b-4a-1)\left(\frac{1}{25}\right)^{b-a-1}.
\end{align*}

{\bf Case 16.} $i=4a+3$ and $j=4b+3$, $0\leq{a}<b\leq{n-1}$. Then
\begin{align*}
z_{i,j}&= \left|
  \begin{array}{cccccccc}
    \frac{4}{5} & -\frac{2}{5} & & & & & &\\
    -\frac{2}{5} & \frac{4}{5} & -\frac{1}{\sqrt{5}} & & & & &\\
    & -\frac{1}{\sqrt{5}} & 1 & -\frac{1}{2} & & & &\\
    & & & \ddots & & & &\\
    & & &  & 1 & -\frac{1}{\sqrt{5}} & &\\
    & & & & -\frac{1}{\sqrt{5}} & \frac{4}{5} & -\frac{2}{5} &\\
    & & & & & -\frac{2}{5} & \frac{4}{5} &  -\frac{1}{\sqrt{5}}\\
    & & & & & & -\frac{1}{\sqrt{5}} & 1\\
  \end{array}
\right|_{(4b-4a-1)\times (4b-4a-1)}\\
&=\frac{2}{25}(4b-4a)\left(\frac{1}{25}\right)^{b-a-1}.
\end{align*}

From Eq. (\ref{sum}), accordingly, we easily obtain the values of $K_{i,j}$ (listed in Eq. (\ref{K_ij})), $0 \leq i<j \leq 3$.

{\bf (1)} Firstly, we calculate $K_{0,0}$ as follows.
\begin{align*}
K_{0,0}&=\sum_{1\leq{a}<b\leq{n}} det\mathcal L_A[4a,4b]=\sum_{1\leq{a}<b\leq{n}} x_{4a-1}y_{4n+1-4b}z_{4a,4b}\\
   &=\sum_{1\leq{a}<b\leq{n}} \frac{1}{6} \left(\frac{1}{25}\right)^{a-1} \cdot
                           \frac{2}{3} \left(\frac{1}{25}\right)^{n-b} \cdot
                           \frac{2}{5}(b-a)\left(\frac{1}{25}\right)^{b-a-1}\\
&=\frac{10}{9}\left(\frac{1}{25}\right)^{n-1} \sum_{1\leq{a}<b\leq{n}} (b-a)
  =\frac{10}{9}\left(\frac{1}{25}\right)^{n-1} \sum_{b=2}^n\sum_{a=1}^{b-1}(b-a)\\
&=\frac{10}{9}\left(\frac{1}{25}\right)^{n-1} \sum_{b=2}^n \frac{b(b-1)}{2}
  =\frac{5}{9}\left(\frac{1}{25}\right)^{n-1} \sum_{b=2}^n (b^2-b)\\
&=\frac{5}{9}\left(\frac{1}{25}\right)^{n-1} \sum_{b=1}^n (b^2-b)
  =\frac{5}{9}\left(\frac{1}{25}\right)^{n-1} \left(\sum_{b=1}^n b^2- \sum_{b=1}^n b\right)\\
&=\frac{5(n^3-n)}{27}\left(\frac{1}{25}\right)^{n-1}.
\end{align*}
Similarly, we obtain the values of $K_{0,2}, K_{0,3}, K_{2,0}, K_{2,2}, K_{2,3}, K_{3,0}, K_{3,2}$, and $K_{3,3}$.

{\bf (2)} Then we calculate $K_{0,1}$.
\begin{align*}
K_{0,1}&=\sum_{1\leq{a}\leq b\leq n} det\mathcal L_A[4a,4b+1]
   = \sum_{1\leq{a}\leq b\leq n-1} det\mathcal L_A[4a,4b+1] + \sum_{1\leq{a}\leq{n}} det\mathcal L_A[4a,4n+1]\\
   &=\sum_{1\leq{a}\leq b\leq{n-1}} \frac{1}{6} \left(\frac{1}{25}\right)^{a-1}\cdot
                           \frac{5}{3} \left(\frac{1}{25}\right)^{n-b}\cdot
                          (4b-4a+1)\left(\frac{1}{25}\right)^{b-a}\\
   &~~~+\sum_{1\leq{a}\leq{n}} \frac{1}{6} \left(\frac{1}{25}\right)^{a-1} \cdot1\cdot
                          (4n-4a+1)\left(\frac{1}{25}\right)^{n-a}\\
&=\frac{5}{18}\left(\frac{1}{25}\right)^{n-1} \sum_{1\leq{a}\leq b\leq{n-1}} (4b-4a+1)
   +\frac{1}{6}\left(\frac{1}{25}\right)^{n-1}\sum_{1\leq{a}\leq{n}}(4n-4a+1)\\
&=\frac{5(4n^3-9n^{2}+5n)}{108}\left(\frac{1}{25}\right)^{n-1}
  +\frac{2n^2-n}{6}\left(\frac{1}{25}\right)^{n-1}\\
&=\frac{20n^3-9n^2+7n}{108}\left(\frac{1}{25}\right)^{n-1}.
\end{align*}
The values of $K_{1,1}, K_{2,1}$, and $K_{3,1}$ can be obtained similarly.

{\bf (3)} Finally, we obtain the values of $K_{1,0}, K_{1,2}$, and $K_{1,3}$ in an almost same way.
For example, $K_{1,0}$ is calculated as follows.
\begin{align*}
K_{1,0}&=\sum_{0\leq{a}<b\leq{n}} det\mathcal L_A[4a+1,4b]
  =\sum_{1\leq{a}<b\leq{n}} det\mathcal L_A[4a+1,4b]
    +\sum_{1\leq{b}\leq{n}} det\mathcal L_A[1,4b]\\
&=\sum_{1\leq{a}<b\leq{n}} \frac{5}{3} \left(\frac{1}{25}\right)^{a} \cdot
     \frac{2}{3} \left(\frac{1}{25}\right)^{n-b} \cdot
     \frac{1}{4}(4b-4a-1)\left(\frac{1}{25}\right)^{b-a-1}\\
&~~~+ \sum_{1\leq{b}\leq{n}} 1\cdot\frac{2}{3} \left(\frac{1}{25}\right)^{n-b} \cdot
                          \frac{1}{4}(4b-1)\left(\frac{1}{25}\right)^{b-1}\\
&=\frac{5}{18}\left(\frac{1}{25}\right)^{n-1} \sum_{1\leq{a}<b\leq{n}} (4b-4a-1)
  +\frac{1}{6}\left(\frac{1}{25}\right)^{n-1}\sum_{1\leq{b}\leq{n}}(4b-1)\\
&=\frac{5(4n^3-3n^{2}-n)}{108}\left(\frac{1}{25}\right)^{n-1}
  +\frac{2n^2+n}{6}\left(\frac{1}{25}\right)^{n-1}\\
&=\frac{20n^3+21n^2+13n}{108} \left(\frac{1}{25}\right)^{n-1}.
\end{align*}


\begin{thebibliography}{99}
\small \setlength{\itemsep}{-.8mm}

\bibitem{b1} H. Wiener, Structural determination of paraffin boiling points, {\em J. Am. Chem. Soc.} {\bf 69} (1947) 17--20.

\bibitem{b2} I. Gutman, Selected properties of the schultz molecular topological index, {\em J. Chem. Inf. Comput. Sci.} {\bf 34} (1994) 1087--1089.

\bibitem{b3} D.J. Klein, M. Randi\'{c}, Resistance distance, {\em J. Math. Chem.} {\bf 12} (1993) 81--95.

\bibitem{b4} H. Chen, F. Zhang, Resistance distance and the normalized Laplacian spectrum, {\em Discrete Appl. Math.} {\bf 155} (2007) 654--661.

\bibitem{b5} H. Zhu, D.J. Klein, I. Lukovits, Extensions of the Wiener number, {\em J. Chem. Inf. Comput. Sci.} {\bf 36} (1996) 420--428.

\bibitem{b6} I. Gutman, B. Mohar, The quasi-Wiener and the Kirchhoff indices coincide, {\em J. Chem. Inf. Comput. Sci.} {\bf 36} (1996) 982--985.

\bibitem{cycle} D.J. Klein, I. Lukovits, I. Gutman, On the definition of the hyper-wiener index for cycle-containing structures, \emph{J. Chem. Inf. Comput. Sci.} \textbf{35} (1995) 50--52.

\bibitem{circulant} H. Zhang, Y. Yang, Resistance distance and Kirchhoffindex in circulant graphs, \emph{Int. J. Quantum Chem.} \textbf{107} (2007) 330--339.

\bibitem{regular} J.B. Liu, X.F. Pan, F.T. Hu, The Laplacian polynomial of graphs derived from regular graphs and applications, {\em Discrete Appl. Math.} {\bf 157} (2014) 18--29.

\bibitem{composite} H. Zhang, Y. Yang, C. Li, Kirchhoffindex of composite graphs, \emph{Discrete Appl. Math.} \textbf{157} (2009) 2918--2927.

\bibitem{complete-mult} R.B. Bapat, M. Karimi, J.B. Liu, Kirchhoff index and degree Kirchhoff index of complete multipartite graphs, \emph{Discrete Appl. Math.} \textbf{232} (2017) 41--49.

\bibitem{flower} N. Faught, M. Kempton, A. Knudson, Resistance distance, Kirchhoff index, and Kemenys constant in flower graphs, {\em J. Math. Chem.} {\bf 86} (2021) 405--427.

\bibitem{strong} Z. Li, Z. Xie, J. Li, Y. Pan, Resistance distance-based graph invariants and spanning trees of graphs derived from the strong prism of a star, {\em Appl. Math. Comput} {\bf 382} (2020) 125335.

\bibitem{Cartesian} J.B. Liu, X.B. Peng, J.J. Gu, W. Lin, The (multiplicative degree-) Kirchhoff index of graphs derived from the Cartesian product of $S_n$ and $K_2$, {\em J. Math.} {\bf 2022} 1670984.

\bibitem{nano} S. Li, W. Sun, S. Wang, Multiplicative degree-Kirchhoff index and number of spanning trees of a      zigzag polyhex nanotube TUHC $[2n; 2]$, {\em Int. J. Quantum Chem.} {\bf 119} (2019) e25969.

\bibitem{LinearPheny} Y. Yang, Computing the Kirchhoff index of linear phenylenes, \emph{J. Combin. Math. Combin. Comput.} \textbf{81} (2012) 199--208.

\bibitem{LinearPheny2} Y. Peng, S. Li, On the Kirchhoff index and the number of spanning trees of linear phenylenes, \emph{MATCH Commun. Math. Comput. Chem.} \textbf{77} (2017) 765-780.

\bibitem{CyclicPheny} L. Ye, On the Kirchhoff index of cyclic phenylenes, \emph{J. Math. Study} \textbf{45} (2012) 233-240.

\bibitem{MobiusPheny} X. Geng, P. Wang, L. Lei, S. Wang, On the Kirchhoff indices and the number of spanning trees of M\"{o}bius phenylenes chain and cylinder phenylenes chain, \emph{Polycycl. Aromat. Comp.}, 2019. DOI: 10.1080/10406638.2019.1693405

\bibitem{CylinderPheny} X. Ma, H. Bian, The normalized Laplacians, degree-Kirchhoff index and the spanning trees of cylinder phenylene chain, \emph{Polycycl. Aromat. Comp.}, 2019. DOI: 10.1080/10406638.2019.1665553

\bibitem{Linear-n-Pheny} S. Li, W. Wei, S. Yu, On normalized Laplacians, multiplicative degree-Kirchhoff indices, and spanning trees of the linear [n] phenylenes and their dicyclobutadieno derivatives, \emph{Int. J. Quantum Chem.} \textbf{119} (2019) e25863.

\bibitem{GeneralPheny} C. Liu, Y. Pan, J. Li, On the Laplacian spectrum and Kirchhoff index of generalized phenylenes, \emph{Polycycl. Aromat. Comp.}, 2019. DOI: 10.1080/10406638.2019.1703765

\bibitem{GeneralPheny2} Z. Zhu, J. B. Liu, The normalized Laplacian, degree-Kirchhoff index and the spanning tree numbers of generalized phenylenes, \emph{Discrete Appl. Math.} \textbf{254} (2019) 256-267.

\bibitem{Periodic} A. Carmona, A. M. Encinas, M. Mitjana, Kirchhoff index of periodic linear chains. \emph{J. Math. Chem.} \textbf{53} (2015) 1195-1206.

\bibitem{Linear4} J. Huang, S. Li, X. Li, The normalized Laplacian, degree-Kirchhoff index and spanning trees of the linear polyomino chains, \emph{Appl. Math. Comput.} \textbf{289} (2016) 324-334.

\bibitem{LinearCrossed4} Y. Pan, C. Liu, J. Li, Kirchhoff indices and numbers of spanning trees of molecular graphs derived from linear crossed polyomino chain, \emph{Polycycl. Aromat. Comp.}, 2020. DOI: 10.1080/10406638.2020.1725898

\bibitem{Linear5} C. He, S. Li, W. Luo, Calculating the normalized Laplacian spectrum and the number of spanning trees of linear pentagonal chains, \emph{J. Comput. Appl. Math.} \textbf{344} (2018) 381-393.

\bibitem{Linear6} Y. Yang , H. Zhang, Kirchhoffindex of linear hexagonal chains, \emph{Int. J. Quantum Chem.} \textbf{108} (2008) 503-512.

\bibitem{Linear6-2} J. Huang, S. Li, L. Sun, The normalized Laplacians, degree-Kirchhoff index and the spanning trees of linear hexagonal chains, \emph{Discrete Appl. Math.} \textbf{207} (2016) 67-79.

\bibitem{LinearCrossed6} Y. Pan, J. Li, Kirchhoff index, multiplicative degree-Kirchhoff index and spanning trees of the linear crossed hexagonal chains, \emph{Int. J. Quantum Chem.} \textbf{118} (2018) e25787.

\bibitem{Mobius6} X. Ma, H. Bian, The normalized Laplacians, degree-Kirchhoff index and the spanning trees of hexagonal M\"{o}bius graphs, \emph{Appl. Math. Comput.} \textbf{355} (2019) 33-46.

\bibitem{Linear8} Q. Zhu, Kirchhoff index, degree-Kirchhoff index and spanning trees of linear octagonal chains, \emph{Australas. J. Comb.} {\bf 153} (2020) 69--87.

\bibitem{LinearCrossed8} J. Zhao, J.B. Liu, S. Hayat, Resistance distance-based graph invariants and the number of spanning trees of linear crossed octagonal graphs, \emph{J. Appl. Math. Comput.} \textbf{63} (2020) 1--27.


\bibitem{Mobius8} J.B. Liu, T. Zhang, Y. Wang, W. Lin, The Kirchhoff index and spanning trees of M\"{o}bius/cylinder octagonal chain, {\em Discrete Appl. Math.} {\bf 307} (2022) 22--31.

\bibitem{Linear8-4} J.B. Liu, J. Zhao, Z. Zhu, On the number of spanning trees and normalized Laplacian of linear octagonal-quadrilateral networks, \emph{Int. J. Quantum Chem.} \textbf{119} (2019) e25971.

\bibitem{Mobius8-4} J.B. Liu, Q. Zheng, S. Hayat, The normalized Laplacians, degree-Kirchhoff index, and the complexity of M\"{o}bius graph of linear octagonal-quadrilateral networks, {\em J. Math.} {\bf 2021} (2021) 1--25.

\bibitem{Decom} J. Huang, S. Li, On the normalized Laplacian spectrum degree-Kirchhoff index and spanning trees of graphs, \emph{Bull. Aust. Math. Soc.} \textbf{91} (2015) 353-367.


\bibitem{TreeNum} N. Biggs, N.L. Biggs, B. Norman. Algebraic graph theory. Cambridge university press, 1993.
%
%
%
%
%
%
%
%
%
%
%
%
%
%
%
%
%
%
%

\end{thebibliography}
\end{document}